\documentclass[11pt]{article}

\usepackage{amsfonts}
\usepackage{amsmath}
\usepackage{upgreek}
\usepackage{amssymb}
\usepackage{amscd}
\usepackage{amsthm}
\usepackage{url}

\newcommand{\ed}{\end{document}}
\newcommand{\tX}{\mathfrak{X}^{\mathrm{sym}}}

\newcommand{\T}{{\tt T}}
\newtheorem{thm}{Theorem} 
\newtheorem{thm*}{Theorem}

\newtheorem{lem}[thm]{Lemma}

\newtheorem{rem*}{Remark}

\numberwithin{equation}{section}

\begin{document}
\title{On the Nonexistence of Skew-symmetric Amorphous Association Schemes}
\author{
Dedicated to the Memory of Robert A. Liebler\\
\\
Jianmin Ma\\ Oxford College of Emory University\\
Oxford, GA 30054, USA\\ Jianmin.Ma@emory.edu}

\date{}
\maketitle




\abstract{
An  association scheme  is amorphous if it has as many fusion schemes as possible. Symmetric amorphous schemes were classified by A. V. Ivanov [A. V. Ivanov,  Amorphous cellular rings II, in  Investigations in algebraic theory of combinatorial objects, pages 39--49. VNIISI, Moscow, Institute for System Studies, 1985] and commutative amorphous schemes were classified by T. Ito, A. Munemasa and M. Yamada  [T. Ito, A. Munemasa and M. Yamada,  Amorphous association schemes over the Galois rings of characteristic 4, European J. Combin., 12(1991), 513--526]. A scheme is called skew-symmetric if the diagonal relation is the only symmetric relation.  We prove the nonexistence of skew-symmetric amorphous schemes with at least 4 classes. We also prove that non-symmetric amorphous schemes are commutative.
} 

\section{Introduction}\label{s:1}

Let ${\mathfrak X}=(X, R=\{R_0, R_1, \ldots , R_d\})$ be a commutative association scheme with 
$d$ classes.  ${\mathfrak X}$
is called amorphous if it has as many fusion schemes as possible. If ${\mathfrak X}$ is symmetric, then it is 
amorphous if and only if every partition of  $R$ containing $\{R_0\}$ gives rise to a fusion scheme.  
 However, if  ${\mathfrak X}$ is non-symmetric,  
 then in order for a partition of $R$ containing  $\{R_0\}$  to give 
rise to a fusion scheme, this partition has to be closed under taking inverse, 
i.e., it is admissible \cite{ItM91}. So, if ${\mathfrak X}$ is non-symmetric, then it is 
amorphous if and only if every admissible partition of  $R$ gives rise to a fusion scheme.

A. V. Ivanov (\cite{Ivanov85}, see also \cite{Gol94})  classified symmetric amorphous 
association schemes with at least three classes: 
 all basic graphs in such a scheme  are strongly regular graphs of Latin square types, or  they 
are all negative Latin square type. Association schemes with two classes are amorphous by definition and
there are many examples in  which the basic graphs are  not either  Latin square types. 
Hence the assumption ``at least three classes'' is essential.
T. Ito, A. Munemasa and M. Yamada \cite{ItM91} classified commutative amorphous  association
schemes  under the assumption $\uptheta + \upphi \ge 3$, where $\uptheta$ is the number of 
pairs of non-symmetric relations, and  
$\upphi$ is the number of non-diagonal symmetric relations (see Section~\ref{s:prelim}). 
The assumption $\uptheta + \upphi \ge 3$ garantees that their symmetrizations have at least three 
classes.  
What about 
association schemes with $\uptheta + \upphi\le 2$? 


This paper  addresses the case $(\uptheta,  \upphi)=(2,0)$. Included in this case are  four-class  
association schemes which have no non-diagonal symmetric relations.
Association schemes with this property will be  referred to as \emph{skew-symmetric}. 
The symmetrizations of skew-symmetric schemes with four classes, as we will see,  indeed have basic
graphs of Latin square type or negative Latin square type. However, such schemes can not exist, due to 
the following theorem. 

\begin{thm} \label{t:main}
There is no  skew-symmetric  amorphous association scheme with 4 classes. 
\end{thm}

Surprisingly, this simple result eleminates the existence of many amorphous association schemes.
\begin{thm}[Main Theorem]\label{t:main2}
There is no skew-symmetric amorphous association scheme with at least 4 classes.
\end{thm}

Theorem~\ref{t:main} answers a question put forward by E. Bannai and S.Y. Song  (\cite[p.395]{Ban93}) regarding the existence of  certain amorphous association schemes with 4  classes.
The proofs of both theorems rely on Theorem \ref{t:pr20}, in which we determine the eigenmatrices of 
skew-symmetric schemes with 4 classes.   

We note that  all association schemes in other cases of $\uptheta +  \upphi  < 3$ are trivially 
amorphous. Let $\mathfrak{X}$ be an association scheme with $\uptheta + \upphi <3$. If
$(\uptheta,  \upphi)=(0,1)$, $\mathfrak{X}$ is  a complete graph.
 If $(\uptheta,  \upphi)=(1,0)$, $\mathfrak{X}$ is a doubly regular tournament. If
 $(\uptheta,  \upphi)=(0,2)$, $\mathfrak{X}$ is equivalent to a pair of complementary strong
 regular graphs. Many chapters of books have been devoted to strong regular 
graphs (e.g. \cite{God93}). If $(\uptheta,  \upphi)=(1,1)$, $\mathfrak{X}$ is a non-symmetric
association scheme with 3 classes, and examples of primitive ones are  not abundant except 
the Liebler-Mena family \cite{Liebler87} and some  examples in \cite{Ma07}
 (see \cite{Jorgensen} and the references there). 
 
It is natural to ask if there exist non-commutative amorphous schemes. We rule out this possibility with 
an algebraic argument.

The general references are  
\cite{BaI84,God93} for association schemes and strongly regular graphs, and 
\cite{Gol94,ItM91,vanDam08} for amorphous association schemes.  In the rest of this paper, all association schemes
are assumed to be commutative unless otherwise stated.

\vspace{12pt}
\textbf{Acknowledgment}: The author would like to thank Professors Robert A. Liebler and Kaishun Wang for many helpful 
discussions and suggestions while preparing this paper. This paper is revised according to the referee's report on
an earlier version, and the author is indebted to the anonymous referee for their valuable comments.   
He is also grateful  to Professor Misha Klin for many encouragements. 

\section{Preliminaries}\label{s:prelim}

Let $X$ be a finite set with cardinality $n \ge 2$  and $R=\{R_0, R_1, \ldots , R_d\}$ be a set 
of binary relations on $X$.  ${\mathfrak X} = (X, R)$ is called an \emph{association scheme with $d$ classes} 
(\emph{a $d$-class association scheme}, or simply, \emph{a scheme}) if the following axioms are satisfied: 
\renewcommand{\theenumi}{\roman{enumi}}
\renewcommand{\labelenumi}{(\theenumi)}
\begin{enumerate}
\item 
$R$ is a partition of $X\times X$ and $R_0 = \{(x,x)\;|\;x \in X\}$ is the diagonal relation.
\item  
For $i=0,1,\dots, d$, the inverse ${R_i^\T} = \{(y, x) \vert (x, y) \in  R_i\}$ of $R_i$ 
is also among the relations:  ${R_i^\T} = R_{i'}$ for some $i'$ ($0 \le  i' \le d$).
\item 
For any triple of $i, j, k = 0, 1, \dots, d$, there exists an integer $p^k_{ij}$ such that for all
$(x,y) \in R_k$,
\[
 |\{z \in  X\;|\; (x, z) \in  R_i, (z, y)\in  R_j\}| = p^k_{ij}.
\] 
\end{enumerate}  
The integers $p^k_{ij}$ are called the \emph{intersection numbers}. 
The integer $k_i = p^0_{i i'} $ is called the \emph{valency} of $R_i$. 
In fact, for any $x\in X$, $k_i = |\{y\in X \vert\; (x,y) \in R_i\}|$.

Furthermore, ${\mathfrak X}$ is called \emph{commutative} if  $ p^k_{ij} =p^k_{ji}$ for all $i, j, k.$

$R_i$ and $R_{i'}$ are called \emph{paired relations}. If $i = i'$, then $R_i$ is called 
\emph{symmetric} or \emph{self-paired}. Let 
$$
\uptheta = |\{ \{i, i'\}| i\ne i', 1\le i \le d\}|, \quad
 \upphi= |\{ i | i=i' , 1\le i \le d\}|.
$$
$\mathfrak{X}$ is called \emph{symmetric} if all relations $R_i$ are symmetric:  $\uptheta =0$. Otherwise, 
$\mathfrak{X}$ is said to be \emph{non-symmetric}. We call $\mathfrak{X}$ 
\emph{skew-symmetric}  if $R_0$ is the only self-paired relation: $\upphi =0$.  

A partition 
$\Lambda_0, \Lambda_1,\ldots, \Lambda_e$ of index set $\{0,1,\dots,d\} $ of $R$
is called \emph{admissible} \cite{ItM91} if $\Lambda_0 =\{0\}$, $\Lambda_i \ne \emptyset$ and 
$\Lambda_i^\T = \Lambda_j$  for some $j\ (1\le i, j\le e)$, where  
$\Lambda^\T =\{\alpha' | \alpha \in \Lambda\}$ is called the inverse of $\Lambda$. 
We may also talk about these properties in terms of 
the relations  when it is convenient, which we did at 
the beginning of the Introduction.

Let
$R_{\Lambda_i} = \cup_{\alpha\in \Lambda_i} R_\alpha$.  If  
$(X, \{R_{\Lambda_i}\}_{i=0}^e)$ is an association scheme, it is called
 a \emph{fusion scheme} of ${\mathfrak X}$.  In particular, the fusion scheme
${\tX} = (X, \{R_0, R_i \cup R_i^\T\}_{i=0}^d)$ is called the \emph{symmetrization}
 of $\mathfrak{X}$. 
 ${\mathfrak X}$ is \emph{amorphous} if every admissible partition 
gives a fusion scheme. Note that if ${\mathfrak X}$ is symmetric, then every partition  containing
$\{0\}$ is admissible by definition. Amorphous  schemes are extremal in the sense they 
have as many fusion schemes as possible.

Let us recall the (first) eigenmatrix $P=(P_{ij})$ of a commutative association scheme
$\mathfrak{ X} = (X, \{R_i\}_{i=0}^d)$.  
Let 
 $A_i$ and $E_i$ $(0\le i\le d)$ be the \emph{adjacency matrices} and \emph{primitive idempotents}
 of  $ \mathfrak{ X}$.   Then the  eigenmatrix $P$  is a square matrix of order $d+1$ defined by 
$$
  A_j = \sum_{i=0}^d  P_{ij} E_i \mbox{ for } j = 0,1,\dots, d.
  $$
 $P$ is characterized by   $A_j E_i = P_{ij}E_i$ for  all 
$i,j =0, 1, \dots d$.  We may index the rows and columns of ${P}$ by 
  $E_i$ and $A_i$ ($i = 0, 1,2$),  respectively.
Moreover, $P_{0i} = k_i $ and $P_{i 0} = 1$ $(0\le i \le d)$. Let $m_i= \mbox{rank } E_i$.   
Then $A_j$ has eigenvalues $P_{0 j} = k_j, P_{1 j},\dots, P_{d j}$ with multiplicities
$m_0=1, m_1,\dots m_d$, respectively. 
The rows and columns of $P$ satisfy the \emph{orthogonality relations}:
\begin{equation}\label{e:orth}
\sum_{i=0}^d \frac{1}{k_i}P_{ji}\overline{P}_{ki} = \frac{n}{m_j}\delta_{jk},\qquad
\sum_{i=0}^d m_iP_{ij}\overline{P}_{ik} = n k_j \delta_{jk},
\end{equation}
where $\overline{x}$ is the complex conjugate of $x$ and $\delta$ is the Kronecker symbol. 

The numbers $m_j$ and  $p_{ij}^{\ell}$  can be calculated from $P$: 
\begin{equation}\label{e:mult}
m_j = \dfrac{n} {
\sum\limits_{i=0}^d \dfrac{1}{k_i} P_{ji} \overline{P}_{ji}
} ,
\end{equation}
\begin{equation}\label{e:int}
p_{ij}^{\ell} = \frac{1}{n k_{\ell}}
      \sum_{h =0}^d m_h P_{hi} P_{hj} \overline{P}_{h\ell }. 
\end{equation}

In the rest of this paper,  the following theorem, referred as the Bannai-Muzychuk criterion for fusion schemes, 
 will be used repeatedly (see \cite{Ban93}, \cite{ItM91}).

\begin{thm}\label{t:BM} 
Let
$\mathfrak{ X} = (X, \{R_i\}_{i=0}^d)$ be a commutative association scheme. Let  
$\{\Lambda_i\}_{i=0}^e$ be an admissible partition of the index set 
$\{0,1,\dots,d\} $. Then  $\{\Lambda_i\}_{i=0}^e$ gives rise to a fusion scheme
$(X, \{R_{\Lambda_i}\}_{i=0}^e)$ if and only if there exists a dual partition  
$\{\Lambda^*_i\}_{i=0}^e$ of $\{0,1,\dots,d\} $ with $\Lambda_0^* =\{0\}$ 
such that each $(\Lambda^*_i, \Lambda_j)$ 
block of the  eigenmatrix $P$ has constant row sum. Moreover, the constant row sum of the 
$(\Lambda^*_i, \Lambda_j)$-block is the $(i,j)$ entry of the eigenmatrix of the fusion scheme. 
\end{thm}

Now we consider two-class association schemes. Let    
${\mathfrak X}=(X,\{R_0, R_1, R_2\})$ be an association scheme.
If ${\mathfrak X}$ is non-symmetric, then $R_2 = {R_1^\T}$.  Its eigenmatrix is 
\begin{equation}\label{e:ch2} 
P = \left[ \begin{array}{ccc} 
1  &  k  & k \\ 1 & \dfrac{-1 + \sqrt{-n}}{2}  & \dfrac{-1 - \sqrt{-n}}{2}\\
1 & \dfrac{-1 - \sqrt{-n}}{2}  & \dfrac{-1 + \sqrt{-n}}{2}
\end{array}\right],
\end{equation}
where $n = |X|  = 2 k+1$. 

Two-class symmetric  schemes are closely related to strongly regular graphs.  
A regular graph $(X,F)$ with  vertex set $X$, edge set $F$ and valency $k$, is called  
\emph{strongly regular} if  
$$
\lambda = |\{z | (x,z)\in F, (z,y)\in F\}|
$$ is constant for all $(x,y)\in F$
 and 
$$
\mu = |\{z | (x,z)\in F, (z,y)\in F\}|
$$
 is constant for all $(x,y)\notin F \; (x\ne y)$. 
The numbers $n=|X|, k, \lambda,\mu$ are the parameters of this graph. 

A strongly regular graph with parameters $(n,k, \lambda, \mu)$ is of (positive) 
\emph{Latin square type} or \emph{negative Latin square
type} if $n=v^2$ (a square) and either (i)
$$ k = g(v-1), \quad \lambda = (g-1)(g-2)+ v -2, \quad \mu = g(g-1),
$$
or (ii)
$$k = g(v+1), \quad \lambda = (g+1)(g+2) - v -2, \quad \mu = g(g+1).
$$
They are denoted  by $L_g(v)$ and $NL_g(v)$, respectively. Graphs with 
$L_g(v)$ parameters can be constructed with $g-2$ mutually orthogonal Latin squares 
of order $v$. Graphs with
$NL_g(v)$ parameters do exist: for example, the Clebsch graph is $NL_1(4)$. 

For a symmetric amorphous scheme  ${\mathfrak X}=(X,\{R_i\}_{i=0}^d)$, 
each graph $ (X, R_i) (i\ne 0)$ is strongly regular.
 It was shown in \cite{Ivanov85, Gol94} that if $d \ge 3$, all graphs 
$ (X, R_i) (i\ne 0)$ are strongly regular graphs of Latin square 
type, or they are all  negative Latin square type. 
The converse is also true (\cite[Theorem 3]{vanDam03}): 
if $(X, R_i)$, $ i=1,\dots, d$ are  strongly regular graphs of all Latin square type or 
 all negative Latin square type  such that $\cup_{i=1}^d R_i = X\times X - R_0 $ and 
$R_i \cap R_j = \emptyset$ ($i\ne j$), then   $(X, \{R_i\}_{i=0}^d)$ 
 is an amorphous association scheme.
 In \cite{ItM91},
T. Ito et al.  classified commutative amorphous association schemes 
with $\uptheta+\upphi \ge 3$ and determined   
 their eigenmatrices and intersection numbers. They also constructed some 
 amorphous schemes on  Galois rings of characteristic 4. 

If $(X,R_1)$ is a strongly regular graph with parameters $(n, k_1, \lambda,\mu)$, 
 the complement $(X, {R_2})$ of $(X,R_1)$ 
is also strongly regular, where ${R_2} = X\times X -R_0- R_1$. 
Furthermore,  $(X,\{R_0,R_1,R_2\})$ is a  symmetric scheme, which
has  the following  eigenmatrix:
\begin{equation} \label{e:srg}
 P = \left[\begin{array}{ccc} 
      1 & k_1 & k_2 \\ 1  & r & t \\ 1 & s & u 
      \end{array}\right]
      \begin{array}{l} 1 \\ m_1 \\ m_2,\end{array}
\end{equation}
where $t = -r -1$, $u = -s-1$. The numbers $k_1, r, s$ are the eigenvalues of the 
adjacency matrix $A_1$ of $R_1$ and  
$r, s$ may be expressed in terms of $n, k_1, \lambda$ and $\mu$. 
Here, we write the multiplicities to the right of $P$. Conversely, a two-class symmetric
scheme gives rise to a pair of complementary strongly regular graphs.

Now we conclude this section with two lemmas that we will need later. In the rest of this paper, we always choose $r \ge 0 > s$. 
\begin{lem} \label{L:char}
Let $\Gamma$ be a strong regular graph with eigenvalues $k, r, s$ and multiplicities 
$1, m_1$, $m_2$. 
\begin{enumerate}
\item [\em{(i)}]
If $m_1 = m_2$ (hence $k=m_1$),  $\Gamma$ is a strong regular graph with parameters $n = 4\mu+1$, 
$k = 2\mu$, $\lambda = \mu -1$. Such a graph is called a conference graph, denoted by $C(n)$. 
\item [\em{(ii)}]
If $k = m_1$, $\Gamma$ is $L_g(v)$ with $v=r-s$, $g=-s$.
\item   [\em{(iii)}]
If $k = m_2$, $\Gamma$ is  $NL_g(v)$ with $v = r-s$, $g = r$. 
 \end{enumerate}
\end{lem} 
 One can prove 
this lemma directly using \cite[exercise 5, p.244]{God93}) or see  Theorem 2.1 of \cite{Mesner67}.
If $v$ is odd, we note that $L_{\frac{1}{2}(v+1)}(v)$ and $NL_{\frac{1}{2}(v-1)}(v)$ have identical
parameters and both agree with $C(n)$ with the argument $\mu = (v^2-1)/4$.

We state the next lemma without a proof since it is straightforward. 
\begin{lem}\label{L:conj}
Let $ \mathfrak{X}$  be a $d$-class association scheme with 
adjacency matrices  $A_i$  and primitive idempotents $E_i$.  
\begin{enumerate}
\item [\em{(i)}]
If $A_i^\T = A_j$, then $P_{\alpha {i}} = \overline{P}_{\alpha{j}}$ $(0\le \alpha \le d)$. So, if
$A_i^\T \ne A_i$, $A_i$ has at least one pair of nonreal eigenvalues that are complex conjugates.  
\item [\em{(ii)}]
If $E_i^\T = E_j$, then $P_{i\alpha} = \overline{P}_{j\alpha}$ $(0\le \alpha \le d)$. So, if 
$E_i^\T = E_i$, then $P_{ik}$ are real for all $k$ $(0\le k \le d)$, and if $E_i^\T \ne E_i$, there are 
distinct $\alpha,\beta$ such that $P_{\alpha i}$ and $P_{\beta i}$ are nonreal complex conjugates.
 \end{enumerate}
\end{lem}
\section{Eigenmatrices}\label{s:table}

Let $\mathfrak{X} = (X, \{R_i\}_{i=0}^4)$ be 
a skew-symmetric association scheme: $(\uptheta,\upphi) = (2,0)$. 
 It is commutative since any association scheme  with at most 4  classes 
is commutative \cite{Hig75}.  Up to a permutation, we may assume $R_4= R_1^\T$ and $R_3= R_2^\T$. Let 
${\tX} =(X,\{R_0, R_1\cup R_4, R_2\cup R_3\})$, 
 the symmetrization of  $\mathfrak{X}$.
 We will determine the eigenmatrix of $\mathfrak{X}$ from that of
${\tX}$. 

Let $\widetilde{A}_i$ and $\widetilde{E}_i$ ($0\le i\le 2$) 
be the adjacency matrices and the primitive idempotents of ${\tX}$, respectively.
Suppose   that  the eigenmatrix $\widetilde{P}$ of  ${\tX}$ 
has form (\ref{e:srg}):
$$ \widetilde{P} = \left[\begin{array}{ccc} 
      1 & k_1 & k_2 \\ 1  & r & t \\ 1 & s & u 
      \end{array}\right] 
\begin{array}{c} 1 \\ m_1 \\ m_2 \end{array}  .
$$
In the rest of this paper,  we  assume that $r \ge 0> s$.

S. Y. Song \cite{Song96} mentioned that up to permutation of rows and columns,  
 a feasible eigenmatrix of $\mathfrak{X}$ can be described as follows: 
\begin{equation} \label{e:p4}
P= \left[\begin{array}{ccccc} 
      1 & k/2  & (n-k-1)/2 & (n-1-k)/2 & k/2\\ 
      1  & \rho  & \tau & \overline{\tau} & \overline{\rho}\\ 
	 1 & \sigma  & \omega & \overline{\omega}& \overline{\sigma} \\
	 1 & \overline{\sigma} & \overline{\omega} & \omega & \sigma\\
     1 & \overline{\rho}  & \overline{\tau} & \tau & \rho 
      \end{array}\right]   
\begin{array}{c} 1 \\ m \\ (n-m-1)/2 \\ (n-m-1)/2\\ m \end{array},
\end{equation}
where  the pair $\rho$ and $\omega$ or the pair $\tau$  and $\sigma$ are nonreal.  He also 
gave a one-sentence explanation. Here, we will prove  Song's observation and determined the entries of $P$.
Let $P$ be the eigenmatrix of $\mathfrak{X}$: 
$$
P = [P_{ij}]_{ 0 \le i, j \le 4}.
$$

Since $A_4 = A_1^\T$, by Lemma \ref{L:conj},  $A_1$ has at least one pair of nonreal eigenvalues. 
 There are two cases to consider:

(1)  $A_1$ has precisely one pair of nonreal  eigenvalues
 $\rho, \overline{\rho}$.  
By Theorem \ref{t:BM},  $\rho + \overline{\rho} = r$ or $s$. 

Suppose $\rho + \overline{\rho} = r$. We may arrange the primitive idempotents  $E_i$ of $\mathfrak{X}$ 
such that
$A_1 E_1 = \rho E_1$,  $A_4 E_4 = \overline{\rho} E_4.$ Hence $E_4 = E_1^\T$. By Lemma \ref{L:conj}, 
$P_{21}, P_{24}, P_{31}$ and $P_{34}$ are all real, and $P_{21}=P_{24}, P_{31}=P_{34}$.  Consider the 
remaining two  primitive idempotents $E_2, E_3 \ne E_0$. Then we have either  $E_3= E_2^\T$ or 
 $E_i = E_i^\T\; (i = 2, 3)$,  and the latter can not occur as we will see.

 Suppose    $E_i = E_i^\T\; (i = 2, 3)$. 
By Lemma \ref{L:conj}, the second and third rows of $P$ have all real entries. 
Since $ \widetilde{A}_1 = A_1 + A_4 $, $P_{21} + P_{24},   P_{31} + P_{34} \in \{r,s\}$ again by Theorem \ref{t:BM}.
We must have $P_{21} + P_{24}\ne    P_{31} + P_{34}. $
 Otherwise, the second and third rows of $P$ are identical, which
 contradicts that $P$ is nonsingular.  We may assume without loss of generality that $P_{21}= P_{24}= r/2$,  
$P_{31} = P_{34} = s/2$.   Since $A_2$ and $A_3$ have a pair of nonreal eigenvalues $\tau$ and
$\overline{\tau}$, $P$ has the following form:
$$
P= \left[\begin{array}{ccccc} 
      1 & n_1 & n_2 & n_2 & n_1\\ 
      1  & \rho  & \tau & \overline{\tau} & \overline{\rho}\\ 
	 1 & r/2  & t/2 & t/2& r/2 \\
	 1 & s/2 & u/2 & u/2 & s/2\\
     1 & \overline{\rho}  & \overline{\tau} & \tau & \rho 
      \end{array}\right]
\begin{array}{c} 1\\  f_1 \\ f_2 \\ f_3 \\ f_1\end{array}, 
$$
where $n_i$ are the valencies of $\mathfrak{X}$ and $2n_i = {k}_i$.
 Now we calculate the multiplicity $f_2$ using (\ref{e:mult}): 
$$
f_2 = \frac{n} { \sum\limits_{i=0}^4 \frac{1}{n_i}P_{2i} \overline{P}_{2i} } = 
\frac{n} { 1 + \frac{2} {n_1} \left(\frac {r}{2}\right)^2  + 
\frac{2} {n_2} \left(\frac {t}{2}\right)^2 } 
= \frac{n} { 1 + \frac{r^2} {k_1}  + 
\frac{t^2} {k_2} }\ ,
$$
which is $m_1$ by (\ref{e:mult}).  Similarly, we can obtain $f_3 = m_2$. So, 
$f_1= 0$, impossible.  

Now we have $E_3= E_2^\T$. 
Since $E_2^\T \ne E_2$ and $P_{21}, P_{41}$ are real, by Lemma \ref{L:conj}, 
$P_{22}, P_{23}$ are nonreal and $P_{22}= \overline{P}_{23}$. Since $A_2^\T= A_3$,
  $P_{32} = \overline{P}_{22} = \overline{P}_{33}$ again by Lemma \ref{L:conj}. 
 Let $\omega = P_{22}$.   So $P$ has the following form:
\begin{equation} \label{e:P1}
P= \left[\begin{array}{ccccc} 
      1 & n_1 & n_2 & n_2 & n_1\\ 
      1  & \rho  & \tau & \overline{\tau} & \overline{\rho}\\ 
	 1 & s/2  & \omega & \overline{\omega}& s/2 \\
	 1 & s/2 & \overline{\omega} & \omega & s/2\\
     1 & \overline{\rho}  & \overline{\tau} & \tau & \rho 
      \end{array}\right],
\end{equation}
where $\rho$ and $\omega$ are both nonreal, $\rho + \overline{\rho} = r,$ 
$\omega+\overline{\omega} = u$ and $\tau + \overline{\tau} =t$.

Suppose $\rho + \overline{\rho} = s$. Replacing $s$ by $r$ in (\ref{e:P1}), we obtain  matrix $P$ for this case, 
in which $\rho$ and $\omega$ are both nonreal, $\rho + \overline{\rho} = s$, 
$\tau + \overline{\tau} = u$ and $\omega + \overline{\omega} = t$. 

(2) $A_1$ has two pairs of nonreal  eigenvalues: $\rho, \overline{\rho}$, and 
$\sigma, \overline{\sigma}$. 

Without loss of generality, we may assume that $\rho + \overline{\rho} = r$,  
$\sigma + \overline{\sigma} =s$, and the first column of $P$ is $(n_2, \rho, 
\sigma, \overline{\sigma}, \overline{\rho})^\T$. So  $P$ has   form (\ref{e:p4}), 
where $\tau$ and $\omega$ can not be both real by Lemma \ref{L:conj}. 

We know from the above analysis that $P$ has the  form asserted in  
(\ref{e:p4}). Let 
$$ 
P= \left[\begin{array}{ccccc} 
      1 & n_1 & n_2 & n_2 & n_1\\ 
      1  & \rho  & \tau & \overline{\tau} & \overline{\rho}\\ 
	 1 & \sigma  & \omega & \overline{\omega}& \overline{\sigma} \\
	 1 & \overline{\sigma} & \overline{\omega} & \omega & \sigma\\
     1 & \overline{\rho}  & \overline{\tau} & \tau & \rho 
      \end{array}\right]
\begin{array}{c} 1 \\ f_1 \\ f_2 \\ f_2 \\ f_1 \end{array}.
$$ 
Now we are ready to  determine $\rho,\omega, \tau$ and $\sigma$. 
 Set 
$$ \rho = \frac{1}{2} (r + \sqrt{-y}), \quad
   \tau = \frac{1}{2} (t + \sqrt{-z}), \quad
   \sigma= \frac{1}{2} (s + \sqrt{-b}), 
$$
where $y, z, b \ge 0$. So,  the first row and column of $P$ are fixed, and hence 
the fourth row and column by Lemma \ref{L:conj}. 
Therefore,  $\omega= \frac{1}{2} (u \pm \sqrt{-c})$ for some $c \ge 0$. 
There are two cases to consider: 

\paragraph*{Case (i)} $\omega = \frac{1}{2} (u + \sqrt{-c})$.
Applying the first orthogonality relation to the first row of $P$, we obtain
$$
1 + \frac{2\rho\overline{\rho}}{n_1} + \frac{2\tau\overline{\tau}}{n_2} = \frac{n}{f_1}.
$$
Substituting $\rho$ and $\tau$ into the above equation, we obtain
\begin{equation} \label{eq:o1}
1 + \dfrac{1}{2 n_1} (r^2 +y) + \dfrac{1}{2n_2} (t^2 + z) = \frac{n}{f_1}. 
\end{equation}
Applying the first orthogonality relation to the first row of $\widetilde{P}$, we obtain
\begin{equation} \label{eq:o2}
1 + \dfrac{r^2}{k_1}  + \dfrac{t^2}{k_2} = \dfrac{n}{m_1}.
\end{equation}
Note that $2n_i= k_i$ and  $2f_i= m_i$. 
By (\ref{eq:o1}) and (\ref{eq:o2}), we have 
\begin{equation} \label{eq:yz}
 \frac{y}{k_1} + \frac{z}{k_2} = \frac{n}{m_1}.
\end{equation}
Similarly,  we can obtain
\begin{equation} \label{eq:bc}
 \frac{b}{k_1} + \frac{c}{k_2} = \frac{n}{m_2}.
\end{equation}

Applying the first orthogonality relation to the second and third rows of $P$ and  $\widetilde{P}$, 
 we obtain
$$
1 + \frac{1}{n_1}(\rho\overline{\sigma} + \overline{\rho}\sigma) + 
      \frac{1}{n_2}(\tau\overline{\omega}+\overline{\tau}\omega) =0, \qquad
1 + \frac{rs}{k_1} + \frac{tu}{k_2}= 0.      
$$
Substituting $\rho, \omega, \tau$, $\omega$ and the second equation into the first equation, we obtain
$$
\frac{\sqrt{by}}{k_1} + \frac{\sqrt{cz}}{k_2} =0.
$$
It follows that $by=0$ and $cz = 0$.  
Since $P$ is nonsingular, we must have either $(y,c)=(0,0)$ or $(z,b)=(0,0)$, but not
both.

Suppose  $(z,b) = (0,0)$. By equations (\ref{eq:yz}) and (\ref{eq:bc}),
\begin{equation}\label{eq:1b}
y = \dfrac{nk_1}{m_1}, \quad c = \dfrac{nk_2}{m_2}.
\end{equation}

Suppose  $(y,c)=(0,0)$. By equations (\ref{eq:yz}) and (\ref{eq:bc}), 
\begin{equation}\label{eq:1a}
 z = \frac{nk_2}{m_1}, \quad  b = \frac{nk_1}{m_2}.
\end{equation}

Note that  in each case above, $\widetilde{P}$ determines $P$ uniquely. 

\paragraph*{Case (ii) } $ \omega = \frac{1}{2} (u - \sqrt{-c})$.
If $c=0$, then $y=0$. This case has been treated in Case (i).  Now, assume $c> 0$. 
Since $\omega $ is nonreal, $ \rho$ is also nonreal: $y > 0$. 
 Note that equations (\ref{eq:yz}) and (\ref{eq:bc}) still hold.  

Now, applying the first orthogonality relation to the second and third rows of  $P$ and simplifying it in 
a similar way,  we can obtain 
\begin{equation}\label{eq:bcyz}
\frac{\sqrt{by}}{k_1} - \frac{\sqrt{cz}}{k_2} =0.
\end{equation}  

Note if $z=0$, then $b=0$. This has been handled previously. Now we assume $z > 0$. Hence,
 $b >0$. Therefore, 
$\rho, \omega, \tau$, and $ \sigma$ are all nonreal, and $b, c, y$, and $z$ satisfy 
equations (\ref{eq:yz}), (\ref{eq:bc}), and (\ref{eq:bcyz}).

\smallskip

We summarize the above discussion in the following theorem. 
\begin{thm} \label{t:pr20}
Let $\mathfrak{X} = (X, \{R_0, R_1, R_2, R_{3}, R_{4}\})$ be a 
skew-symmetric association scheme with $R_1 = R_4^\T$ and $R_2 = R_3^\T$.
Let $\widetilde{P}$ be  the eigenmatrix of the symmetrization  ${\tX}$: 
$$ \widetilde{P} = \left[\begin{array}{ccc} 
      1 & k_1 & k_2 \\ 1  & r & t \\ 1 & s & u 
      \end{array}\right]
\begin{array}{c} 1 \\ m_1 \\ m_2\end{array},
$$
where $r \ge 0 > s$. 
Then the eigenmatrix of $\mathfrak{X}$ has the following form:
$$ 
P= \left[\begin{array}{ccccc} 
      1 & k_1/2  & k_2/2 & k_2/2 & k_1/2\\ 
      1  & \rho  & \tau & \overline{\tau} & \overline{\rho}\\ 
	 1 & \sigma  & \omega & \overline{\omega}& \overline{\sigma} \\
	 1 & \overline{\sigma} & \overline{\omega} & \omega & \sigma\\
     1 & \overline{\rho}  & \overline{\tau} & \tau & \rho 
      \end{array}\right].
$$
 $\rho,\omega, \tau$ and $\sigma$ take values in one of three cases:
\begin{enumerate}
\item [\em{(i)}] 
 $ 
  \sigma = \dfrac{s}{2}, \quad \tau = \dfrac{t}{2}, \quad
 \rho =  \dfrac{1}{2}\left( r + \sqrt{-\dfrac{nk_1}{m_1}} \right), \quad 
\omega =\dfrac{1}{2}\left( u + \sqrt{-\dfrac{nk_2}{m_2}} \right).
$

\item [\em{(ii)}] 
$ \rho = \dfrac{r}2, \quad \omega = \dfrac{u}{2}, \quad 
   \sigma= \dfrac{1}{2} \left(s + \sqrt{- \dfrac{nk_1}{m_2}}\right), \quad  
   \tau = \dfrac{1}{2} \left(t + \sqrt{- \dfrac{nk_2}{m_1}}\right). 
$
\item  [\em{(iii)}]
$ \rho = \frac{1}{2} \left(r + \sqrt{-y}\right), \quad
   \tau = \frac{1}{2} \left(t + \sqrt{-z}\right), \quad
   \sigma= \frac{1}{2} \left(s + \sqrt{-b}\right), \quad
   \omega = \frac{1}{2} \left(u - \sqrt{-c}\right),
$\\[5pt]
where $b,c,y$, and $z$ are all positive and  satisfy the following equations: 
$$
\frac{y}{k_1} + \frac{z}{k_2} = \frac{n}{m_1}, \quad
 \frac{b}{k_1} + \frac{c}{k_2} = \frac{n}{m_2}, \quad
 \frac{\sqrt{by}}{k_1} - \frac{\sqrt{cz}}{k_2} =0.
 $$
\end{enumerate}
\end{thm}

\section{Proof of the Main Theorem}\label{s:main}

In this section, we prove Theorems \ref{t:main} and \ref{t:main2}. 
Let $\mathfrak{X} = (X, \{R_0, R_1, R_2, R_{3}, R_{4}\})$ be a skew-symmetric amorphous  scheme,
whose eigenmatrix $P$ is given in Theorem \ref{t:pr20}. 

Suppose that  $P$  is given by Theorem~\ref{t:pr20}(i). 
Since $\mathfrak{X}$ is amorphous,  $R_0, R_1\cup R_2, R_{3}\cup R_{4}$ gives rise to a
skew-symmetric association  scheme, which has eigenmatrix (\ref{e:ch2}).
By Theorem \ref{t:BM},
 $$
\rho + \tau  =\dfrac{1}{2}\left( r + t+ \sqrt{-\dfrac{nk_1}{m_1}}\right)
=\dfrac{-1 + \sqrt{-n}}{2}.
$$ So $m_1=k_1$. By Lemma \ref{L:char}(ii),  
$(X,  R_1\cup R_4)$ is $L_g(v)$ with $v = r-s$, $ g =-s$.

The $i$-th intersection matrix $B_i$ is a square matrix of order $d+1$ whose 
$(j,\ell)$ entry is $p_{ij}^{\ell}$. Using (\ref{e:int}), we can obtain
$$
B_1 =\left[ \begin {array}{ccccc} 
0 & 1 & 0 & 0 & 0 \\
0 & \dfrac{\lambda + r}{4}&
\dfrac{k_1 (k_1 -\lambda -1 -t)}{4k_2}&
\dfrac{k_1 (k_1 -\lambda -1 -t)}{4k_2}&
\dfrac{\lambda - 3r}{4}\\\noalign{\medskip}
0 & \dfrac{k_1 -\lambda -1 +t}{4} &
\dfrac{k_1 -\mu + s}{4}&
\dfrac{k_1 -\mu - s}{4}&
\dfrac{k_1 -\lambda-1 - t}{4}\\\noalign{\medskip}
0 & \dfrac{k_1 -\lambda -1 +t}{4}&
\dfrac{k_1 -\mu - s}{4}&
\dfrac{k_1 -\mu + s}{4}&
\dfrac{k_1 -\lambda-1 - t}{4}\\\noalign{\medskip}
\dfrac{k_1}{2} & \dfrac{\lambda + r}{4}&
\dfrac{k_1 (k_1 -\lambda -1 +t)}{4k_2}&
\dfrac{k_1 (k_1 -\lambda -1 +t)}{4k_2}&
\dfrac{\lambda + r}{4}
\end {array} \right] .
$$ 
From $p_{12}^1 - p_{12}^4 = {2t}/{4}$ and
$p_{12}^2 - p_{12}^3 = {2s}/{4}$
 we readily 
deduce that $t$ and $s$ are   even integers. Since $r+t +1 =0$, $r$ is odd. On the other hand, 
  $p_{11}^1 = \frac{\lambda + r} {4} =   \frac{s(s + 2)  + 2r }{4}$. So $r$ is even because
$s(s+2)$ is divisible by 4, a contradiction.

Suppose that  $P$  is given by
Theorem~\ref{t:pr20}(ii). In a similarly way, we can deduce  
 $k_1 =m_2$ and hence $(X,  R_1\cup R_4)$ is $NL_g(v)$ with $v = r-s$, $ g= r$. The first 
intersection matrix for this case
 can be obtained from the above $B_1$  by interchanging $r$ and $s$, and $t$ and $u$. 
We can readily deduce that   $r$ and $u$ are  even
and  $s$ is odd. Since
$p_{11}^1 =  \frac{\lambda + s}{4} = \frac{r(r+2) + 2s }{4}$, $s$ is even,
 a contradiction.

Suppose  that $P$  is given by 
Theorem \ref{t:pr20}({iii}). In this case, $\rho, \omega, \tau, \sigma$ are all nonreal. 
since $\mathfrak{X}$ is amorphous, 
 $R_0, R_1\cup R_2, R_3\cup R_4$ determines a non-symmetric association scheme, 
which has eigenmatrix (\ref{e:ch2}). So,  
$$ \rho + \tau = \dfrac{1}{2} (r+t + \sqrt{-y} + \sqrt{-z}) = \dfrac{-1 + \sqrt{-n}}{2},$$ 
and hence $\sqrt{-y} + \sqrt{-z} =  \sqrt{- n}$. 
Similarly, $(X, \{R_0, R_1\cup R_3, R_2\cup R_4\})$ is non-symmetric association scheme and thus
$\sqrt{-y} - \sqrt{-z} = \pm  \sqrt{- n}$. These  equations 
imply either $y =0$ or $z =0$, a contradiction. 
This completes the proof of Theorem \ref{t:main}. \hfill $\Box$

Now we prove the main theorem.  Suppose that $\mathfrak{X}$ is a 
skew-symmetric amorphous scheme with more than 4  classes. 
$\mathfrak{X}$ has a skew-symmetric fusion scheme with 4 classes,  
 which can not exist by Theorem \ref{t:main} because this fusion scheme is again amorphous. 
This completes the proof.

\section {Concluding Remarks}

As we mentioned in the Introduction that there does not exist any non-commutative amorphous scheme. 
We now give a short proof of this result. By the main theorem, we may assume $\upphi \ge 1$. 
Since association schemes with at most 4 classes are commutative, we may assume $\uptheta+\upphi \ge 3$.
We first treat the minimal cases $(\uptheta,\upphi)=(1,3)$ or $
 (2, 1)$ and  the  general case will then follow.  
Suppose that $\mathfrak{X}$ is  a non-symmetric amorphous 
scheme with $(\uptheta,\upphi)=(1,3)$ or $
 (2, 1)$. If $\mathfrak{X}$ is non-commutative, then the adjacency algebra generated by 
 the adjacency matrices of $\mathfrak{X}$ over the complex numbers $\mathbb{C}$
 is non-commutative of dimension 6. It is semisimple  and thus 
is isomorphic to direct sum of full matrix algebras of degree 1,1 and 2:
$$
\mathbb{C} \oplus \mathbb{C} \oplus M_2(\mathbb{C}).
$$

$\mathfrak{X}$  has a 4-class fusion scheme $\mathfrak{F}$, which is commutative.  So
the  adjacency algebra of $\mathfrak{F}$ is commutative of dimension 5. On the other hand,   
$M_2(\mathbb{C})$ can not have commutative subalgebras of dimension 3 and thus the adjacency 
algebra of $\mathfrak{X}$ has no commutative subalgebras of dimension 5, which is a contradiction. 
Therefore, $\mathfrak{X}$ is commutative. (In fact, we have proved that a non-commutative scheme 
with 5 classes can not have a 4-class fusion scheme.)

Let $\mathfrak{X}$ be an amorphous non-symmetric scheme with $\uptheta + \upphi > 4, \uptheta \ge 1$.
Any two adjacency matrices of $\mathfrak{X}$ are among the adjacency matrices of some fusion scheme with
$(\uptheta,\upphi)=(1,3)$ or $
 (2, 1)$. So the adjacency matrices of $\mathfrak{X}$ commute pairwisely and thus $\mathfrak{X}$ is commutative.

We note that the minimal cases can also be handled by a careful analysis of their intersection numbers, 
which shows that the intersection numbers in each  case coincide with those of 
certain commutative amorphous scheme. 
In fact, using the notation in \cite{ItM91}, 
amorphous schemes with $(\uptheta,\upphi)=(1,3)$  belong to $L_{g_1;g_2,g_3, g_4}(v)$ or 
$NL_{g_1;g_2,g_3, g_4}(v)$, and 
 amorphous schemes with $(\uptheta,\upphi)=(2,1)$ belong to $L_{g_1,g_2;g_3}(v)$ or 
$NL_{g_1,g_2;g_3}(v)$.

We conclude with some remarks:
\begin{itemize}

\item  E. R. van Dam and M. Muzychuk \cite{vanDam08} gave an excellent survey of symmetric 
amorphous association schemes. Among many results, they gave all known constructions and enumeration of such schemes with
 vertices up to 49 vertices. However, there has not been  much work done with the non-symmetric counterpart except 
\cite{ItM91}.  

\item 
In light of Theorem~\ref{t:pr20}, it is interesting to study skew-symmetric schemes with 4 classes. 
We are working on it in another paper. 
\item
In the literature (e.g. \cite{Gol94},\cite{Ivanov85},\cite{vanDam08}), 
an association scheme is call \emph{amorphic}  if every partition of  $R$  containing $\{R_0\}$  
gives rise to a fusion scheme. 
 The notion of  admissible partition  was introduced and  the term amorphous was used by T. Ito, et al. in \cite{ItM91}. For symmetric schemes, the two notions are equivalent. 
 Any association scheme with two classes is trivially amorphic 
by definition. It is easy to see that any amorphic scheme with at least three classes is symmetric. 

\item 
 I.N. Ponomarenko an A.R. Barghi \cite{Ponomarenko07} recently introduced amorphic
$C$-algebras by axiomatizing the property that each partition of a standard basis leads to 
a fusion algebra. Just like association schemes, each amorphic $C$-algebra of 
dimension $\ge 4$ is symmetric. 
They  showed that each amorphic 
$C$-algebra  is determined up to isomorphism by the multiset of its degrees 
(valencies in the case of association scheme) and an additional integer
$\epsilon = \pm{1}$ (reflecting the positive or negative Latin square type). 
Since our focus here is 
non-symmetric association schemes, our work has little overlap with that of
 \cite{Ponomarenko07} .

\end{itemize}


\end{document}